\magnification=\magstep1
\input amstex
\documentstyle{amsppt}
\widestnumber \key {GHMR}
%\NoRunningHeads
\NoBlackBoxes

%\Macros-- Definitions
\let\text=\hbox

     %Tate-Shafarevich Group

% Common words in a formula
\def\and{\text{ and }}

\def\ff{\text{ if }}   % \if is already in the macros of the computer. Can %not define it again

%\def\oo{\text{ or }}  % \or  is already in the macros of the computer. Can %not define it again

  %Tamagama number

%Common math words

  %inflation

\def\mod{ \text{ mod }}
\def\ord{\text{ ord}}
  % Reduced

\def\tr{\text{ tr}}

%Common math notations

\def\IM{\hbox{ Im}}

\def\Gal{\text{ Gal}}

%NAMES

 \def\A{\Bbb A}

\def\O{\Cal O}

\def\lr{\longrightarrow}

\parskip 4pt plus 1.5pt\baselineskip 26pt plus 0pt minus 0pt\lineskip 6pt plus 1pt \lineskiplimit 5pt

\baselineskip=24 true pt \hsize=5.75 true in \hoffset=.565 true in
\vsize=8.65 true in \voffset=.125 true in

\parskip 4pt plus 1.5pt\baselineskip 26pt plus 0pt minus 0pt\lineskip 6pt plus 1pt \lineskiplimit 5pt

\baselineskip=24 true pt \hsize=5.75 true in \hoffset=.565 true in
\vsize=8.65 true in \voffset=.125 true in

\def\odd{\hbox{odd}}

\topmatter
\title
On CM abelian varieties over imaginary quadratic fields
\endtitle
%\rightheadtext{Nonvanishing of Central Hecke L-value}
\author
Tonghai Yang
\endauthor
\address
Department  of Mathematics University of Wisconsin Madison, WI
53706 USA
\endaddress
\email thyang\@math.wisc.edu
\endemail
\abstract
 In this paper, we associate canonically to every imaginary quadratic field
 $K=\Bbb Q(\sqrt{-D})$  one or two isogenous
 classes of CM abelian varieties over $K$,
 depending on whether $D$ is odd  or even ($D \ne 4$). These abelian
varieties are characterized as of smallest dimension and smallest
conductor, and such that the abelian varieties themselves descend
to $\Bbb Q$. When $D$ is odd or divisible by 8, they are the
`canonical' ones first studied by Gross and Rohrlich. We prove
that these abelian varieties have the striking property that the
vanishing order of their $L$-function at the center is dictated by
the root number of the associated Hecke character. We also prove
that the smallest dimension of a CM abelian variety over $K$ is
exactly the ideal class number of $K$ and classify when a CM
abelian variety over $K$ has  the smallest dimension.
\endabstract

\thanks
Partially supported by an AMS Centennial Fellowship and  a  NSF
grant DMS-0070476
\endthanks
\endtopmatter

\document

\subheading{0. Introduction} The paper is motivated by two basic
questions related to   an  imaginary quadratic field $K=\Bbb
Q(\sqrt{-D})$ with fundamental discriminant $-D$ and abelian
varieties over $K$ with complex multiplications (CM). It is
well-known that there is a CM elliptic curve over $K$ if and only
if $K$ has ideal class number $1$. what is then the smallest
dimension of a CM abelian variety over $K$ in general? We prove
that the smallest dimension is exactly the ideal class number $h$
of $K$ (Theorem 3.1).  We also classify the CM abelian varieties
over $K$ of dimension $h$ in terms of its associated algebraic
Hecke character (Theorems 3.4 and 3.5) in section 3. It turns out
that it only depends on the restriction of the Hecke character on
the principal ideals. To be more precise, let $(A, i)$ be a CM
abelian variety over $K$ of CM type $(T, \Phi)$, and let $\chi$ be
the  associated  algebraic Hecke character of $K$ of conductor
$\frak f$, then
$$
\chi(\alpha \O_K) = \epsilon(\alpha) \alpha \tag{0.1}
$$
for some odd character $\epsilon: (\O_K/\frak f)^* \lr \Bbb C^*$.
We will prove in section 3

\proclaim{Theorem 0.1} Let the notation be as above, and let $h$
be the ideal class number of $K$. Then

(1) (Theorem 3.1) \quad  one has $h | \dim A$.

(2) (Theorem 3.5) \quad  The equality $\dim A =h$  holds if and
only if one of the following holds.

 (a) \quad $\epsilon$ is quadratic,

  (b) \quad $3|D$ and $\epsilon$ is of order $6$, there is $\alpha
  \in K^*$ such that $3 \alpha^2$ is prime to $\frak f$ and
   $\epsilon(3 \alpha^2) = -1$, or

   (c) \quad $4|D$ and $\epsilon$ is of order $4$ or $12$, there
   is $\alpha
  \in K^*$ such that $2 \alpha^2$ is prime to $\frak f$ and
   $\epsilon(2 \alpha^2) = \pm i$.

   In such a case, $A$ is a scalar restriction of a CM elliptic
   curve over the Hilbert class field $H$ of $K$ if and only if
   $\hbox{Im}(\epsilon) \subset \O_K^*$ (Proposition 2.2).
   \endproclaim

   In view of the theorem, it is nature to ask whether and how
   one can
   associate `canonically' a `nice' CM abelian variety over $K$
   of dimension $h$ to an imaginary quadratic field $K$.  When $D$ is odd or $8||D$, this can be done by
   means of the scalar restriction of certain `canonical' CM
    elliptic  curves over $H$, according to  Gross ([Gr]) when $D$ is an odd prime
      and Rohrlich ([Ro2-4]) in general.
   These CM abelian varieties  descend to varieties over $\Bbb Q$ and
   have bad reductions exactly at
   $p|D$.
    Indeed, according to Rohrlich ([Ro2-4]), a
   Hecke characters  $\chi$ of $K$ of conductor $\frak f$
   is called canonical if it
   satisfies the following three conditions:

   (0.2) \quad $\chi(\bar\frak a) = \overline{\chi(\frak a)}$ for all
   ideals $\frak a$ of $K$ prime to $\frak f$;

    (0.3) \quad The character $\epsilon$ in $(0.1)$  is quadratic;

    (0.4) \quad The conductor $\frak f$ is divisible only primes
       ramified in $K/\Bbb Q$.

If $\chi$ is a canonical Hecke character of $K$, so is $\chi \phi$
for every ideal class character  $\phi$ of $K$  or $\chi^\sigma$
for every $\sigma \in \Gal(\bar Q/K)$ with the same $\epsilon$.
They form the same family and gives rise to the same CM abelian
variety $A$ over $K$ up to isogeny (different characters
correspond to different embeddings $i$ and different CM types).
The condition $(0.2)$ implies that $A$ descends to an abelian
variety over $\Bbb Q$, and $(0.3)$ implies that $A$ is a scalar
restriction of some CM elliptic curve over $H$. $(0.2)$ and
$(0.4)$ imply that $A$ has bad reduction exactly at primes
dividing $D\O_K$. One of the most striking properties of the
canonical Hecke characters (or `canonical CM abelian varieties')
is
$$
\ord_{s=1} L(s, \chi)= \frac{1- W(\chi)}2 \tag{0.5}
$$
where $W(\chi)$ is the root number of $\chi$.

  However, a canonical CM abelian variety over  $K$ exists if and only
  if  $D$ is odd or is divisible by $8$, according to Rohrlich
  ([Ro2]). When $D$ is odd, it is unique up to isogeny, and the
  associated $\epsilon$ is given by
  $(\O_K/\sqrt{-D})^*\cong (\Bbb Z/D)^* \lr \{\pm 1\}$. The
  associated canonical Hecke characters have root number
  $(\frac{2}D)$.
  When $8||D$, there are two canonical CM abelian varieties over
  $K$ up to isogeny, the associated $\epsilon$ are given by
  $$
  \epsilon=\epsilon_2 \epsilon^0: (\O_K/\frak p_2^5)^* \times
  (\O_K/\sqrt{-D_1})^* \lr \{\pm 1\}.
  $$
  Here $\frak p_2$ is the prime ideal of $K$ above $2$, and
  $\epsilon_2$ is a nontrivial quadratic character  of
  $(\O_K/\frak p_2^5)^*$ (two choices), and
  $\epsilon^0=(\frac{}{D_1})$ is the quadratic character of
  $(\O_K/\sqrt{-D_1})^* \cong \Bbb Z/D_1)^*$ and $D_1=D/8$.
  The associated canonical Hecke characters have the root number
  $W(\chi) = \epsilon( 1+\sqrt{-D/4})$. What about $4||D$?
  We will prove in section 4 the following theorem.

   \proclaim{Theorem 0.2}  Let $K =\Bbb Q(\sqrt{-D})$ be a
   quadratic imaginary quadratic field with fundamental discriminant $-D $, and let
   $$
   n(D)=\cases 1 &\ff  D \hbox{ is odd or } D=4,
   \\
     2 &\ff D >4 \hbox{ is even}.
     \endcases
     $$
   Then   there are exactly $n(D)$
   CM abelian varieties $A$ over $K$ of dimension $h$, up to isogeny,
    such that
    the abelian variety itself (not the complex multiplications) descends to
    an abelian variety over $\Bbb Q$ and has the
   smallest   possible  conductor.  Moreover, then $D$ is
   odd or $8||D$, they are the canonical CM abelian varieties
   discussed above.
   \endproclaim
 We call  the CM abelian varieties defined in Theorem 0.2
    the {\it{simplest}} CM abelian varieties over $K$. In sections 5,
    we compute the root numbers of
   corresponding {\it{simplest}} Hecke characters and prove in sections 6 and 7
    that the amazing formula $(0.5)$  holds for all simplest Hecke
    characters of $K$,  which we record  here as

 \proclaim{Theorem 0.3} Let $\chi$ be a simplest Hecke character
 of $K =\Bbb Q(\sqrt{-D})$. Then $(0.5)$ holds.
 \endproclaim

    This theorem, combining with the
    Gross-Zagier, and a deep result of Rubin and/or  Kolyvagin, implies that the simplest CM
    abelian varieties $A$ over $K$ has always finite
    Shafarevich-Tate group over $K$ and has the Mordell-Weil rank
    $0$ or $2h$  over $K$ depending on whether the root number of
    the associated simplest Hecke characters have root number $+1$
    or $-1$. The exceptional case $D=4$ is due to the fact that $i
    \in \Bbb Q(\sqrt{-4})$.

       For
   simplicity, we sometimes  miss out the CM type of CM abelian
   varieties in this paper. This does no harm since the arithmetic or L-function
   of a CM abelian variety does not depend on the choices of the CM types, see Remark 1.2.
   In section 1,  We  review the basic relation between CM abelian
   varieties and the algebraic Hecke characters.
   In section 2, we  study the relation between the
   scalar restriction of CM elliptic curves over the Hilbert class
   field of $K$ and $h$-dimensional CM abelian varieties over $K$.
   In particular, we prove (Proposition 2.2) that a CM abelian
   variety over $K$ of dimension $h$ is a scalar restriction of a
   CM elliptic curve over $H$ if and only if its Hecke character
   has values in $K$ on principal ideals.

{\bf Acknowledgement} This work was inspired by David Rohrlich's
work on canonical Hecke characters and my joint work with Stephen
D. Miller on the same subject. The author thank both of them for
their  inspirations. The author thanks Dick Gross, Stephen Kudla,
Ken Ribet,  and David Rohrlich for stimulus discussions.  Part of
the work was done  while the author visited the Department of
Mathematics at HKUST for 3 weeks in 2002. The author thanks
Jian-Shu Li and the department for their hospitality and financial
support.

 \subheading{1. CM abelian varieties
and algebraic Hecke characters}

Let $K$ be a number field. Then a type of $K$ is a formal sum
$\Phi=\sum_{\sigma: K \hookrightarrow \Bbb C} n_\sigma \sigma$
with $n_\sigma \in \Bbb Z$. If $L$ is a finite extension of $K$,
one can extend $\Phi$ to a type $\Phi_L$ of $L$ via
$$
\Phi_L=\sum_{\sigma: L \hookrightarrow \Bbb C} n_\sigma \sigma,
\quad n_\sigma =n_{\sigma|_K}.
$$
A type $\Phi$ is called simple if it can not be extended from a
proper subfield. When $K$ is Galois over $\Bbb Q$, embeddings of
$K$ into $\Bbb C$ are elements of the Galois group $\Gal(K/\Bbb
Q)$. In this case, one defines its reflex type
$$
\Phi'=\sum_{\sigma: K \hookrightarrow \Bbb C} n_\sigma
\sigma^{-1}.
$$
In general, given a type $(K, \Phi)$, its reflex type $(K',
\Phi')$ is defined as follows. Firstly,  $K'$ is the subfield of
$\Bbb C$  generated by all $\alpha^\Phi=\prod
\sigma(\alpha)^{n_{\sigma}}$, $\alpha \in K^*$. Secondly, let $L$
be a finite Galois extension of $\Bbb Q$ containing both $K$ and
$K'$, extend $\Phi$ to $\Phi_L$. It is then a standard fact that
there is a unique type $\Phi'$ of $K'$ such that its extension to
$L$ is
$$
(\Phi')_L=(\Phi_L)'.
$$
This $\Phi'$ is independent of the choice of $L$ and is called the
reflex type of $\Phi$.

An algebraic Hecke character of $K$ of infinite type $\Phi$ and
modulus $\frak f$ (an integral ideal) is a group homomorphism
$$
\chi: I(\frak f) \longrightarrow \Bbb C^*
$$
such that for every $\alpha \equiv 1 \mod^* \frak f$
$$
\chi(\alpha \O_K) =\alpha^\Phi.
$$
Here $I(\frak f)$ denotes the group of fractional ideals of $K$
prime to $\frak f$, and $\alpha \equiv 1\mod^* \frak f$ means
$$
\ord_v (\alpha-1) \ge \ord_v \frak f \quad \ff \ord_v \frak f >0,
$$
and $\sigma (\alpha) >0$ for every real embedding $\sigma$ of $K$.
The Dirichlet unit theorem implies that there is an algebraic
Hecke character of $K$ of infinite type $\Phi$ if and only if
$w(\Phi)=n_\sigma + n_{\rho \sigma}$ is independent of the choice
of $\sigma$, where $\rho$ is the complex conjugation of $\Bbb C$.
In such a case,  $\Phi$ is called a Serre type of weight
$w(\Phi)$. Notice that the subfield $\Bbb Q(\chi)$ of $\Bbb C$
generated by $\chi(\frak a)$, $\frak a \in I(\frak f)$, is a
number field containing the reflex field $K'$ of $(K, \Phi)$. We
say $\chi$ has values in $T$ if $\Bbb Q(\chi) \subset T$.

  When $K$ is a CM number field, i.e., a quadratic totally
  imaginary extension of a totally real number field, a  type
  $\Phi$ of $K$ is a CM type if $n_\sigma \ge 0$ and
  $n_\sigma + n_{\rho\sigma}=1$ for every complex embedding
  $\sigma$ of $K$. In this case, $\Phi$ is often identified with
  the set of embeddings $\{\sigma: n_\sigma=1\}$. In general, any
  extension of a CM type just defined is also called a CM type.

  An abelian variety $A$ defined over a subfield $L$ of $\Bbb C$
  is said to be a CM abelian variety over $L$ if there is
 a number field $K$ of degree $[K:\Bbb Q]=2 \dim A=2 d$ together
 with  an embedding
  $$
  i: K\hookrightarrow \hbox{End}_L^0A=\hbox{End}_L A \otimes \Bbb
  Q.
  $$
  In such a case, $K$ acts on the differentials $\Omega_{A/\Bbb
  C}$ diagonally via $d$ embeddings $\Phi=\{\phi_1, \cdots,
  \phi_d\}$: there is a basis $\omega_i$ for $\Omega_{A/\Bbb
  C}$ such that
 for every $\alpha \in i^{-1}(\hbox{End}_L
  A)$
  $$i(\alpha)^* (\omega_i)=\phi_i(\alpha) \omega_i.
  $$
  We will identify $\Phi$ with the formal sum $\sum_{\phi\in \Phi}
  \phi$. We usually call $(A, i)$ is of CM type $(K, \Phi)$. It is
  a fact ([Sh2, Theorem 1, Page 40]) that the two seemingly different definitions
  of CM types are the same. We remark that the CM type of $(A, i)$
  depends on $i$, $A$ can have different CM types if one allows
  $i$ change. For example, let $A=E^2$ be the square of a CM
  elliptic curve $E$ by a quadratic field $K$. Then
  $\hbox{End}^0A=M_2(K)$, and any quadratic field extension of $K$,
  embedded into $M_2(K)$, gives rise to CM type of $A$.

   The following theorem summarizes the
  basic relation between CM abelian varieties and algebraic Hecke
  characters.

  \proclaim{Theorem 1.1} Let $K \subset \Bbb C$ be a number
  field,  and let $(T, \Phi)$ be a CM type and $(T', \Phi')$ be
  its reflex type. Assume $T' \subset K$.

  (1) (Shimura-Taniyama)\quad If $(A, i)$ is a CM abelian variety over $K$ of CM
  type $(T, \Phi)$, then there is a (unique) algebraic Hecke
  character $\chi$ of $K$ of infinite type $\Phi'_K$ such that
  $i(\chi(\frak p))$ reduces to the Frobenius endomorphism  of $A$
  modulo $\frak p$ for every prime ideal $\frak p$ where $A$ has good reduction.
  In particular,
  $$
  L(s, A/K)=\prod_{\sigma: T \hookrightarrow \Bbb C} L(s,
  \chi^\sigma).
  $$
  Here $\chi^\sigma= \sigma\circ \chi$.  We call $\chi$ the associated (algebraic) Hecke
  character of $(A, i)$ or simply $A$.

  (2)(Casselman) \quad Conversely, if $\chi$ is an algebraic Hecke
  character of $K$ of infinite type $\Phi'_K$, valued in $T$, then
  there is a CM abelian variety $(A, i)$ over $K$ of CM type $(T,
  \Phi)$, unique up to isogeny, such that the associated Hecke
  character of $(A, i)$ is $\chi$. We call $(A, i)$ or simply $A$
  an associated abelian variety of $\chi$.

  (3)\quad Let $(A, i)$ be a CM abelian variety over $K$ of CM
  type $(T, \Phi)$ and let $\chi$ be the associated Hecke
  character of $K$. Then the following are equivalent.

(a) \quad
  The CM abelian variety $A$ is simple over $K$.

  (b) \quad  $T=\Bbb Q(\chi)$.

%  $(T, \Phi)$ is called primitive if there is a simple abelian variety of CM type
 % $(T, \Phi)$.
  \endproclaim
  \demo{Proof} For (1), see for example [Sh2, Theorems 19.8 and 19.11]. For
  (2), see for example [Sh2, Theorem 21.4].  If $A$ is not simple, then $A,
  $ is isogenous to $B^r$ for some simple abelian variety $B$ over $K$, and
  $B$ is a CM abelian variety of some CM type $(T_1, \Phi_1)$. In
  this case, $\chi_A=j\circ\chi_B$, where $j$ is an embedding of
  $T_1$ to $T$. So $\Bbb Q(\chi_A)\subset T_1 \neq T$. Conversely, if
  $T_1=\Bbb Q(\chi_A)\neq T$, let $(B, i_B)$ be a CM abelian variety
  over $K$ of CM type $(T_1, \Phi_1)$ associated to $\chi_A$,
  where $\Phi_1=\Phi|_{T_1}$. Then $(B^r, i_B')$ is of CM type
  $(T, \Phi)$ for some $i_B'$ where $r =[T:T_1]$. So $(A, i_A)$
  and $(B^r, i_B')$ have the same CM type and the same Hecke character,
  and thus are $K$-isogenous to each other.  $A$ is thus not simple over $K$. This proves
  $(3)$
  \enddemo

\remark{Remark 1.2} If $(A, i)$ is a CM abelian variety of CM type
$(T, \Phi)$, with algebraic Hecke character $\chi$, then $(A,
i\circ \sigma^{-1}) $ is of type $(T^\sigma, \Phi \circ
\sigma^{-1})$, with algebraic Hecke character $\chi^\sigma$. When
$A$ is simple over $K$, $(T^\sigma, \Phi\circ \sigma^{-1})$ and
 $i\circ \sigma^{-1}$ are all possible CM types and embeddings.
 In particular, the arithmetic and L-function of $A$ over $K$ do
 not depend on particular choices of the CM types or the
 embeddings. the abelian variety $A$ is in correspondence with
 the family of algebraic Hecke characters  $\{\chi^\sigma:\,
 \sigma \in \hbox{Aut}(\Bbb C /\Bbb Q) \}$. If one wants to leave
 the infinite type $\Phi_K'$ of $\chi$ unchanged, one needs to
 require $\sigma \in \hbox{Aut}(\Bbb C/K'')$ where $K''$ is the
 reflex field of $(K, \Phi_K')$.
 \endremark
 %We remark in passing that any algebraic Hecke character $\chi$ of
 %$K$ can be written as the product of a Hecke character of $K$ of
 %finite order and powers of Hecke characters associated to CM
 %abelian varieties over $K$.

\subheading{2. CM abelian varieties over a quadratic field and
Scalar restriction of CM elliptic curves}

 From now on, Let $K=\Bbb Q(\sqrt{-D})$ be an imaginary quadratic
 field with fundamental discriminant $-D$, let $H$ be the Hilbert
 class field of $K$, and let $h=[H:K]$ be the ideal class number
 of $K$. Let $E$ be a CM elliptic curve over $H$ of CM type $(K,
 id)$, where $id$ is the prefixed complex embedding of $K$ such
 that $\hbox{Im}(\sqrt{-D}) >0$. Let $\chi_E$ be the
 algebraic Hecke character of $H$ associated to $E$,  with values in $K$. It is
 interesting to study the scalar restriction
 $B=\hbox{Res}_{H/K} E$, which is an abelian variety over $K$ such
 that $B(R)=E(H\otimes_K R)$ for every $K$-algebra $R$. In order for $B$
 to be a CM abelian variety over $K$, $E$ has to be a $K$-curve, in the sense that
 every Galois conjugate $E^\sigma$ for $\sigma \in \hbox{Gal}(H/K)$ is $H$-isogenous to $E$.
 But this condition is not sufficient in general.  We refer
 to [Gr] and   [Na] for results  and references in this direction. In terms
 of Hecke characters, we have

 \proclaim{Proposition 2.1} Let the notation be as above. Then
 $B=\hbox{Res}_{H/K} E$ is a CM abelian variety over $K$ if and
 only if there is an algebraic Hecke character $\chi$ of $K$ such
 that $\chi_E = \chi\circ N_{H/K}$.
 \endproclaim
 \demo{Proof} If $(B, i_B)$ is CM of type $(T, \Phi)$, then
 $\Phi=\{\sigma: T \hookrightarrow \Bbb C:\, \sigma|_K=\hbox{ id
 }\}$. Let $\chi_B$ be the associated algebraic Hecke character of
 $K$ associated to $B$. We have to verify $\chi_E=\chi_B\circ
 N_{H/K}$. They have the same type. So it suffices to prove that
 for a prime ideal $\frak P$ of $H$ where $E$ has good reduction,
 one has
 $$
 \chi_E(\frak P)=\chi_B(N_{H/K}(\frak P)) =\chi_B(\frak p)^f
 $$
 where $\frak p$ is the prime ideal of $K$ below $\frak P$ and
 $f =[\O_H/\frak P: \O_K/\frak p]$. By [Gr, Lemma 15.1.6],
 $$\hbox{End}_K B =\sum_{\sigma \in \Gal(H/K)}
 \hbox{Hom}_H(E^\sigma, E) \sigma.
 $$
 Let $\widetilde{Fr(\sigma)}$ be the Frobenuius map from
 $\widetilde{E^\sigma}$ to $\tilde E$, then
 $\widetilde{Fr}=\sum  \widetilde{Fr(\sigma)} \sigma $ is the
 Frobenuius of $\tilde B$, and thus it is equal to $\chi_B(\frak
 p)$ acting on $\tilde B$. On the other hand, $\tilde\sigma^f=1$,
 $\widetilde{Fr(\sigma)}^f$ is the Frobenuius of $\tilde E$,
 independent of the choice of $\sigma$, and it is thus equal to
 $\chi_E(\frak P)$ acting on $\tilde E$. So we have
 $\chi_E(\frak P)= \chi_B(\frak p)^f$ as desired.

 Conversely, if $\chi_A=\chi\circ N_{H/K}$ for some algebraic
 Hecke character $\chi$ of $K$, then $\chi(\alpha \O_K) \in K^*$
 since every principal ideal of $K$ is a  norm  from $H$. Let
 $T(\chi=\Bbb Q(\chi)$ and $\Phi(\chi)$ be the set of complex  embeddings
 of $T$ extending the fixed embedding of $K$. Let  $A(\chi)$ be
 a CM abelian variety over $K$ of CM type of $(T, \Phi)$
 associated to $\chi$. For every $\sigma \in \Phi$, $\chi^\sigma
 \chi^{-1}$ is trivial  on principal ideals and thus an ideal
 class character. So
 $$
 \{ \chi^\sigma:\, \sigma \in \Phi\} \subset \{\chi\phi: \phi
 \hbox{ is an ideal class character of } K \}.
 $$
 It will be proved in Theorem 3.1 that $\# \Phi \ge h$. So the
 above relation is in fact an equality. This implies
 $L(s, \chi_E)=\prod_{\sigma \in \Phi} L(s, \chi^\sigma)$, and thus $B$ is isogenous to
 $A(\chi)$.
 \enddemo

  According to [Sh2, Theorem 19.13], the condition
  $\chi_E=\chi\circ N_{H/K}$ is also equivalent to that
  $K(E_{tor})$ is the maximal abelian extension of $K$.

\proclaim{Proposition 2.2} Let $(A, i)$ be a simple CM abelian
variety over $K$ of CM type $(T, \Phi)$. Let $\chi_A$ be the
associated  algebraic Hecke character of $K$ with values in $T$.
Then the following are equivalent.

(1) \quad  There is a CM elliptic curve $E$ over $H$ such that $A$
is isogenous to the scalar restriction of $E$ from $H$ to $K$.

(2) \quad  One has  $\chi_A(\alpha \O_K) \in K^*$ for every
principal ideal of $K$ prime to the conductor of $\chi_A$.
\endproclaim
\demo{Proof} $(1)\Rightarrow (2)$. The assumption $(1)$ asserts
that $\chi_H=\chi_A\circ N_{H/K}$. By the global class field
theory, every ideal $\alpha \O_K$ is a norm of some ideal $\frak
A$ of $H$. So $\chi_A(\alpha \O_K) =\chi_E(\frak A) \in K^*$.

$(2) \Rightarrow (1)$. Let $\chi=\chi_A \circ N_{H/K}$. By the
global class field theory again, $N_{H/K}(\frak A)$ is a principal
ideal of $K$ for every ideal $\frak A$ of $H$. So
$$\chi(\frak A) = \chi_A(N_{H/K}(\frak A)) \in K^*,$$ and thus
$\chi$ is an algebraic character of $H$ with values in $K$.
Therefore there is a CM elliptic curve $E$ over $H$ whose
associated algebraic Hecke character is $\chi$. One has
$$
L(s, E/H) =L(s, \chi)L(s, \bar\chi)
          =\prod_\phi L(s, \chi_A \phi)L(s, \overline{\chi_A \phi}).
$$
Here the product runs over all ideal class characters $\phi$ of
$K$. On the other hand, since $A$ is simple, one has $T=K(\chi_A)$
is generated by $\chi_A$ over $K$. So one has by assumption $(2)$
$$
\align
 \{ \chi_A^{\sigma}: \sigma \in \Phi \}
 &=\{\chi_A^{\sigma}:\sigma \in \hbox{Aut}(\Bbb C/K)\}
\\
&=\{ \chi_A \phi: \phi \hbox{ is an ideal class character of }
K\}.
\endalign
$$
So $L(s, A/K)=L(s, E/H)=L(s, \hbox{Res}_{H/K}E)$. So $A$ is
isogenous to $\hbox{Res}_{H/K}E$, proving (1).
\enddemo

\subheading{3. CM abelian varieties over a quadratic imaginary
field of smallest dimension}

 Let $K=\Bbb Q(\sqrt{-D})$ be again a quadratic imaginary field of
 fundamental discriminant $-D$. Let $(A, i)$ be a CM abelian
 variety over $K$ of type $(T, \Phi)$, and let
  $\chi=\chi_A$ be the associated algebraic Hecke character of K of
 conductor $\frak f$. Then there is a character $\epsilon$ of
$(\O/\frak f)^*$ such that
$$
\chi(\alpha \O ) = \epsilon(\alpha) \alpha \tag{3.1}
$$
for every principal ideal of $K$ prime to $\frak f$.  Obviously
$\epsilon(\alpha) =\alpha^{-1}$ for every unit $\alpha$ in K. In
particular $\epsilon(-1)=-1$.

   Let  $K(\chi)$ be the subfield of $\Bbb C$ generated by
   all $\chi(\frak a)$ over $K$, $\frak a \in I(\frak f)$. Then
   $K(\chi) \subset T$,  and $K(\chi) =T$ if and only if $A$ is
   simple.

\proclaim{Theorem 3.1} Let $(A, i)$ be a CM abelian
 variety over $K$ of type $(T, \Phi)$. Then $h| \dim A$.
 \endproclaim

   When $\dim A =1$, this gives the well-known fact that in order
   to have a CM elliptic curve over a quadratic imaginary field,
   the field has to have ideal class number $1$. The simple reason
    is the fact that the $j$-invariant generates the
   Hilbert class field.

\demo{Proof} We may assume that $A$ is simple and thus
$T=K(\chi)$.
 For any number field F, let $\mu(F)$ be the group of
roots of unity in F, and let $w_F =\# \mu(F)$. Choose an integer
$n>0$ such that $w_T|n$ and $\chi(\frak a^n) \in K^*$ for every
$\frak a \in I(\frak f)$. Then
$$
\theta: I(\frak f)/ P(\frak f) \lr K^*/K^{*n}, \quad
       [\frak a]   \mapsto \chi(\frak a^n) K^{*n}
$$
is a well-defined map and is injective. Indeed, for any principal
ideal $\alpha \O \in P(\frak f)$, one has
$$
\chi(\alpha \O)^n = \epsilon(\alpha)^n \alpha^n =\alpha^n \in
K^{*n},
$$
which implies that $\theta$ is well-defined. On the other hand, if
$\chi(\frak a^n) =\alpha^n$ for some $\alpha \in K$. Then $\frak
a^n$ and $\alpha^n $ generates the same ideal in K, so are $\frak
a$ and $\alpha$. Therefore $\theta$ is injective! Clearly, the
image of $\theta$ is in the kernel of the natural map
$$
K^*/K^{*n} \lr T^*/T^{*n}.
$$
So one has  $h|[T:K]$ by [Ro1, Proposition 1] .
\enddemo

 The rest of this section is to determine all CM abelian varieties
 over $K$ of the smallest possible dimension $h$. For this
 purpose, we assume $D >4$ and fix an odd character
 $$
 \epsilon: (O_K/\frak f)^* \longrightarrow   \Bbb C^*, \quad \epsilon(-1)=-1.
 \tag{3.2}
 $$
 Let $\chi$ be an algebraic Hecke character of $K$ satisfying $(3.1)$,
 and let $T=K(\chi)$, and
 $$
 \Phi=\{\sigma: T\hookrightarrow \Bbb C: \, \sigma|_K=id  \}.
 $$
 Let $A=A(\chi)$ be the associated simple CM abelian variety over
 $K$ of type $(T, \Phi)$.

   Let $L=K(\epsilon)$ be the subfield of $T$ generated by
   $\epsilon$ over $K$.
 Let
$$
H(\frak f)=\{ \frak a \in I(\frak f) : \frak a^2 \hbox{ is
principal }\}.
$$
Then $H(\frak f)/P(\frak f)$ is the genus ideal class group of K
and has order $2^r$, where $r+1$ is the number of prime factors of
D. Let $T_g$ be the subfield of T generated by $\chi(\frak a)$,
$\frak a \in H(\frak a)$. We have
$$
K \subset L \subset T_g \subset T. \tag{3.3}
$$

 The same proof as in  [Ro1, Theorem 2] yields

\proclaim{Proposition 3.2}  (1) \quad One has
$$
\align &\{\chi^{\sigma}: \sigma \in \Gal(\bar\Bbb Q/T_g)\}
\\
&=\{\chi^{\sigma}:  \sigma \hbox{ is a complex embedding of T
which is trivial    on } T_g \}
\\
&=\{ \chi \phi: \phi \hbox{ is an ideal class character  trivial
on the genus subgroup} \}.
\endalign
$$

(2) \quad $[T:K]=2^{-r} h [T_g:K]$ depends only on the effect of
$\chi$ on ideals whose square is principal, where $r+1$ is the
number of prime factors of D.

(3) \quad $ \mu_T \subset T_g.
%\tag{H4}
$
\endproclaim
\demo{Proof} (1)\quad The first identity is clear. For a complex
embedding $\sigma$ of T fixing $T_g$, $\phi=\chi^{\sigma}
\chi^{-1}$ is obviously an ideal class character  trivial on the
genus subgroup. So it suffices to prove $\#I(\frak f)/H(\frak f)
\le [T:T_g]$ for the second equality. Let $F=T_g(\mu_T)$, and
$$
\theta: I(\frak f) \lr    F^*/F^{*n}, \frak a \mapsto \chi(\frak
a^n) F^{*n}.
$$
Here $n$  again satisfies $\#\mu_T | n$ and $\chi(\frak a)^n \in
K^*$ for every ideal $\frak a$ of $K$.  Then its image is in the
kernel of the natural map $ F^*/F^{*n} \lr T^*/T^{*n}$, and thus
has order dividing $[T:F]$ by [Ro1, Proposition 1]. Next, define
$$
\psi: \ker \theta \lr K^*/K^{*n}, \frak a\mapsto \chi(\frak a^n)
K^{*n}.
$$
Then its image is in the kernel of the natural map $K^*/K^{*n} \lr
F^*/F^{*n}$.  Since $F$ is an abelian extension of K, [Ro1,
Proposition 2] implies that for any $\frak a \in \ker \theta$
$$
\chi(\frak a^n)^2 \in K^{*n}.
$$
 So $\frak a^2$ is principal, and thus $\ker \theta \subset H(\frak f)$. This proves
$$
\# I(\frak f) /H(\frak f) \le \IM( \theta) \le [T:F] \le [T:T_g].
\tag{3.4}
$$
This proves the second identity of (1), and that the inequalities
in (3.4) are all equalities. In particular,
$$
[T:T_g]=\# I(\frak f) /H(\frak f) =2^{-r} h,
$$
which proves claim (2). One has also $F= T_g$, which gives claim
(3).
\enddemo

\proclaim{Lemma  3.3}  Let, for any prime $p$,
$$
\mu_{p}(\epsilon) =\{ x \in \hbox{Im}(\epsilon): \, x^{p^n}=1
\hbox{ for some } n
 \ge 0 \} .
\tag{3.4} $$
 Then there is a root of unity $c_p \in \mu_{2}(\epsilon)$ for
 each prime number $p|D$ such that
$ T_g=L(\sqrt{c_p p}, p|D)$. Moreover, if $\xi$ is a generator of
$\mu_2(\epsilon)$, then $ T_g \subset L(\sqrt\xi, \sqrt p, p|D)$.
%\hbox{ and } \quad
%\prod_{p|D} c_p \in -1 \cdot L^{*2}.
\endproclaim
\demo{Proof} Let $\frak p_l$ be the prime ideal of $K$ above $l$,
for every $l|D$,  then $\frak p_l$ has order 2 in the ideal class
group of $K$ and generate the genus class group of $K$. For each
$l|D$, choose  $\alpha_p \in K^*$ so that $\frak a_l = \alpha_l
\frak p_l$ is relatively prime to $\frak f$.  Then $\frak a_l$
generate $H(\frak f)/P(\frak f)$. Now
$$
\frak a_l^2= \alpha_l^2 l\O \quad \hbox{ and } \quad \chi(\frak
a_l)=\pm \sqrt{\epsilon(\alpha_l^2 l ) l} \alpha_l \in T_g. $$
 So
$$
\sqrt{\epsilon(\alpha_l^2 l ) l} \in T_g.
\tag{3.5} $$
 Obviously,
$\epsilon(\alpha_l^2 l )$ can be replaced by a root of unity $c_l
\in \mu_2(\epsilon)$. Since $T_g$ is generated by $\chi(\frak
a_l)$ over $L$, we have
$$
L(\sqrt{c_l l}, l|D) = T_g.
$$
Since $\sqrt{c_l} \in L(\sqrt\xi)$, one has thus $T_g \subset
L(\sqrt\xi, \sqrt l, l|D)$.  This proves the lemma.
\enddemo

\proclaim{Theorem 3.4} Let $\epsilon$ be an odd character given by
$(3.2)$.  Let $\chi$ be a Hecke character of K satisfying (3.1),
and let $T=K(\chi)$.

(1) \quad If $[T:K]=h$, then
$$
\IM(\epsilon) \subset
  \cases
    \mu_2 &\ff (6, D) =1,
\\
    \mu_4 &\ff 3 \nmid D,
\\
    \mu_6 &\ff 2 \nmid D,
\\
    \mu_{12} &\ff 6|D.
\endcases
\tag{3.6}
$$

(2)  \quad Conversely, assume that $\epsilon$ satisfies (3.6),
then $[T:K]=h$ or $2 h$, with $[T:K]=h$ if and only if  $\sqrt \xi
\notin T_g$. Here $\xi $ is a generator of $\mu_2(\epsilon)$.

(3) \quad Whether $[T:K]=h$ (equivalently $\dim A(\chi) =h$) or
not depends only on $\epsilon$, i.e., the restriction of $\chi$ on
principal ideals.
\endproclaim

  %Because of Theorem 3.4,  we will call a CM abelian variety $A$ over $K$ of type
  % $\epsilon$ if its associated Hecke character $\chi$ is of
  % type $\epsilon$, i.e., satisfying $(3.1)$.

\demo{Proof} By Proposition 3.2, $[T:K]=h$ if and only if
$[T_g:K]=2^r$.   Let $H_g=K(\sqrt{p^*}, p|D)$ be the genus class
field of K, where
$$
p^*=\cases
   (-1)^{\frac{p-1}2} p &\ff p \ne 2,
  \\
   \frac{-D}{\prod_{p|D, p \ne 2} p^*} &\ff p =2.
  \endcases
%\tag{3.6}
$$
Notice that $[H_g:K]=2^r$, where $r+1$ is the number of prime
factors of D. Since
$$
[L(\sqrt\xi, \sqrt p, p|D): T_g]=1 \hbox{ or } 2
$$
by Lemma 3.3, and
$$
[K(\sqrt\xi ,\sqrt p, p|D):H_g]\ge 2,
$$
one has
$$
\align [T_g:K] &\ge \frac{1}2 [L(\sqrt\xi, \sqrt p, p|D):K]
\\
  &=  2^{r-1}  [L(\sqrt\xi, \sqrt p, p|D):K(\sqrt\xi, \sqrt p, p|D)]
                  [K(\sqrt\xi, \sqrt p, p|D):H_g]
\\
 &\ge 2^r [L(\sqrt\xi, \sqrt p, p|D):K(\sqrt\xi, \sqrt p, p|D)].
\endalign
$$
So $[T_g:K] \ge 2^r$, and $[T_g:K]=2^r$ if and only if the
following three conditions hold.
$$
\align
 &\sqrt \xi \notin T_g,
\tag{3.7}
\\
 &L(\sqrt\xi, \sqrt{ p}, p|D)= K(\sqrt\xi, \sqrt{ p}, p|D),
\tag{3.8}
\\
 &[K(\sqrt\xi, \sqrt{\pm p}, p|D):K(\sqrt{p^*}, p|D)]=2.
\tag{3.9}
\endalign
$$
It is easy to see from (3.9) that $\xi$ has at most order 4, which
occurs only when $4|D$. $(3.8)$ implies then that $L(\sqrt\xi,
\sqrt{ p}, p|D) \subset K(\sqrt{-1}, \sqrt 2 , \sqrt p, p|D)$
contains at most $\mu_{24}$. So $\mu_T=\mu_{T_g} \subset \mu_{12}$
by $(3.7)$. A slightly more careful inspection shows $(3.6)$.

 Now assume that $(3.6)$ is satisfied. Then  $L \subset H_g$ except for two cases: $D=8$ when $\epsilon$ has order $4$, and
 $D=24$ when $\epsilon$ has order $12$. So except for the two
 cases one has
$$
[L(\sqrt\xi, \sqrt p, p|D) : H_g] =[H_g(\sqrt\xi): H_g] =2.
$$
Thus $[T_g:K]=2^r $ or $2^{r+1}$. Moreover $[T_g:K]=2^r$ if and
only if
 $\sqrt\xi \notin T_g$. This proves (2) except for the two exceptional cases. In the exceptional cases,
 $T_p=L$ has degree $2h$ over $K$, and $\sqrt\xi \in T_g$, the claim still holds.  Finally,
  since $c_l=\epsilon(\alpha_l^2 l)$ up to a
 square in $L^*$, $T_g$ is determined by $\epsilon$, thanks to
 Lemma 3.3. Now (3) follows from (2).
\enddemo

\proclaim{Theorem 3.5} Let the notation be as in Theorem 3.4, and
assume that $\epsilon$ satisfies $(3.6)$.  Let $A(\chi)$ be the CM
abelian variety over $K$ associated to $\chi$.

(1) \quad If $\epsilon^2 =1$, then $\dim A(\chi)=h$, and $A(\chi)$
is isogenous to the scalar restriction of a CM elliptic curve $E$
over $H$, the Hilbert class field of $K$.

(2) \quad If $\epsilon$ has order $4$, then $\dim A(\chi) =h$ if
and only if $4|D$, $D\ne 8$, and  for some (and every) element
$\alpha_2 \in K^*$ such that $2 \alpha_2^2$ is prime to $\frak f$,
one has $\epsilon(2 \alpha_2^2)$ is of order $4$.

(3) \quad If $\epsilon$ has order $6$, then $\dim A(\chi) =h$ if
and only if  $3|D$, and for some (and every) element $\alpha_3 \in
K^*$ such that $3 \alpha_3^2$ is prime to $\frak f$ and
$\epsilon(3 \alpha_3^2)$ is of order $2$ or $6$.

(4) \quad If $\epsilon$ has order $12$, then $\dim A(\chi) =h$ if
and only if  $12|D$, and for some (and every)   element $\alpha_2,
\in K^*$ such that $2 \alpha_2^2$ prime to $\frak f$ and
%and  $3 \alpha_3^2$ are prime to $\frak f$ and
$ \epsilon(2 \alpha_2^2)$ is of order $4$ or $12$.
%\quad \epsilon(3 \alpha_3^2) =\pm 1.
\endproclaim

   Clearly, in cases (2)-(4) with $D >4$, $A(\chi)$ is {\it{not}} a scalar
   restriction of any elliptic curve when its dimension is $h$.
  % In the exceptional cases $D=3,4,8$, $K$ has ideal class number
  % $1$, there is nothing to say about it. When $D=24$,
  % $K=\Bbb Q(\sqrt{-24})$ has ideal class number  $2$. In this case, $\dim
  % A(\chi)=2$ if $\epsilon$ has order $2$, $4$ or $6$, and $\dim
  % A(\chi)=4$ when $\epsilon$ has order $12$.

\demo{Proof of Theorem 3.5} (1) follows from Proposition 2.2. For
(2), notice first that $\mu_2(\epsilon)$ is generated by $i$ and
that $4|D$ is necessary by Theorem 3.4(1). Assume thus $4|D$.  If
$\epsilon(2 \alpha_2^2)=\pm 1$, then $\sqrt{\pm 2} =
\sqrt{\epsilon(2 \alpha_2^2)2} \in T_g$ and thus $\zeta_8
=\frac{1}{\sqrt 2}(1+i)\in T_g$. This implies $\dim
A(\chi)=[T:K]=2h$ by Theorem 3.4(2). If $\epsilon(2
\alpha_2^2)=\pm i$, then $\sqrt{\epsilon(2 \alpha_2^2)2}=\sqrt{\pm
2i} =\pm(1\pm i)$, and so
$$
T_g= \Bbb(i,\sqrt{-D/4},  \sqrt{\epsilon_l(l \alpha_l^2) l}, 2\ne
l|D).
$$
Here $\alpha_l \in K^*$ are such that $\frak a_l=\frak p_l
\alpha_l$ is prime to $\frak f$ as in the proof of Lemma 3.3.
Since
$$
\prod_{l|\frac{D}4} \frak a_l = \frac{D}4 \prod_{l|\frac{D}4}
\alpha_l $$ is principal, one has
$$
\prod_{l|\frac{D}4}\sqrt{\epsilon_l(l \alpha_l^2) l}
 =\pm \prod_l \alpha_l^{-1} \chi(\prod_{l|\frac{D}4} \frak a_l)
 \in L^*.
 \tag{3.10}
 $$
This implies that $\sqrt{\epsilon_l(l \alpha_l^2) l}, 2\ne l|D$
have a relation over $L$, and thus $[T_g:K]=2^r$. So $[T:K]=h$.

 For (3), notice first that $3|D$ is necessary and $\mu_2(\epsilon)=\{\pm
 1\}$. Notice also $L=K(\epsilon)=K(\sqrt{-3})$. If $\epsilon(3
 \alpha_3^2) =1$ (up to a cubic root of unity), then $\sqrt 3 \in
 T_g$ by Lemma 3.3, and thus $i=\sqrt{-1} \in T_g$. This implies
 $\dim A(\chi) =2h$ by Theorem 3.4(2). If so $\sqrt{3}$ can
 not be in $T_g$. If $\epsilon(3\alpha_3^2) =-1$ (up to a cubic root of unity),
then $\sqrt{\epsilon(3\alpha_3^2) 3}=\sqrt{-3} \in L$. Then same
argument as in (2) (in particular $(3.10)$) shows $[T_g:K]=2^r$
and thus $\dim A(\chi) =h$.

 For $(4)$, there is $\beta \in K^*$ prime to $D$ such that $\epsilon(\beta)$ and
 $\epsilon(\beta^2)$ has order $3$. Let $\alpha_2 \in K^*$ be such that
 $2\alpha_2^2$ is prime to $D$. Then $\epsilon(2 \alpha_2^2)$ has order $4$
 if and only if $\epsilon(2 \alpha_2^2 \beta^2)$ has order $12$. The rest  is similar
 to $(2)$ and left to the reader.
\enddemo

\proclaim{Corollary 3.6} When $(6, D)=1$, every CM abelian variety
over $K$ of dimension $h$ is $K$-isogenous to a scalar restriction
of a CM elliptic curve over $H$.
\endproclaim

\example{Example 3.7} Assume  that $4|D$ and $D >4$. Let $p \equiv
1 \mod 4$ be a prime number split or ramified in $K=\Bbb
Q(\sqrt{-D})$. Let $\frak p$ be a prime ideal of $K$ above $p$.
Let $\epsilon_p$ be a surjective  character
$$
\epsilon_p: (\O_K/\frak p)^* \cong (\Bbb Z/p)^* \longrightarrow
\mu_4
$$
so that a given generator $g$ of $(\Bbb Z/p)^*$ maps to $i$. Then
$$
\epsilon_p(-1)=-1 \Leftrightarrow p\equiv 5 \mod 8 \Leftrightarrow
\epsilon_p(2) =\pm i.
$$
Let $f=\prod_{p} p $ be the product of finitely many such prime
 numbers $p \equiv  1 \mod 4$, and let $\frak f$ be an integral  ideal of $K$
 whose norm is  $f$.
Let $\epsilon=\prod_p \epsilon_p$ and assume that
$\epsilon(-1)=-1$, i.e., odd number of $p's$ are congruent $5$
modulo $8$. Then $\epsilon(2)=\pm i$ by the above equivalence. So
any CM abelian variety over $K$  of type $\epsilon$ has dimension
$h$ by Theorem 3.5(2).  On the other hand, let $q \equiv 3 \mod 4$
be a prime split or ramified in $K$ and let $\frak q$ be a prime
above $q$. Let
$$
\epsilon_q': (\O_K/\frak q)^* \cong (\Bbb Z/q)^* \lr \{\pm 1 \}
$$
be such that $\epsilon_q'(n) = (\frac{n}q)$ for an integer $n$
prime to $q$. Then  $\epsilon_q'(-1)=-1$. Now let $f'=f q$ with
$f$ being as above but with even number of primes $p \equiv 5 \mod
8$. Let $\epsilon' =\epsilon \epsilon_q'$, then $\epsilon'(-1)=-1$
is odd and $\epsilon'(2) =\pm -1$.  So any CM abelian variety over
$K$ of type $\epsilon$ has dimension $2h$.
    \endexample

\example{Example 3.8}  Similarly, assume that  $3 |D$, and $D >3$.
Let  $p \equiv 1 \mod 3$ be an odd prime number split or ramified
in $K$ and  let $\epsilon_p$ be a surjective  character
$$
\epsilon_p: (\O_K/\frak p)^* \cong (\Bbb Z/p)^* \lr \mu_6.
$$
Then
$$
\epsilon_p(-1)=-1\Leftrightarrow p\equiv 7 \mod 12\Leftrightarrow
3 \hbox{ is a square modulo } p . $$
 Let $f=\prod_p p$ be a
product of odd primes $p \equiv 1 \mod 3$ split or ramified in $K$
and let $\frak f$ be an ideal of $K$ whose norm is $f$. Let
$\epsilon=\prod_p \epsilon_p$ and assume that $\epsilon(-1)=-1$,
i.e., there are odd number of prime divisors of $f$ satisfying $p
\equiv 7 \mod 12$. Then $\epsilon(3)$ is of order $2$ or $6$. So
any CM abelian variety over $K$ of type $\epsilon$ is of dimension
$h$. On the other hand, if we let $q \equiv -1 \mod 12$ and
$\epsilon_q'$ be as in Example 3.7. Let $f$ and $\epsilon$ be as
above but with $\epsilon(-1)=1$, then $\epsilon'=\epsilon
\epsilon_q'$ is odd and any CM abelian variety over $K$ of type
$\epsilon'$ is of dimension $2h$.
 \endexample

\example{Example 3.9} Similarly, assume $ 12|D$. Let $p \equiv 1
\mod 12$ be a prime number split or ramified in $K$. Let
$\epsilon_p$ be a surjective  character
$$
\epsilon_p: (\O_K/\frak p)^* \cong (\Bbb Z/p)^* \lr \mu_{12}.
$$
Define $f$ and $\epsilon$ the same way as in Example 3.7. If
$\epsilon(-1)=-1$, then any CM abelian variety over $K$ of type
$\epsilon$ is of dimension $h$. if $\epsilon(-1) =1$, let $q$ and
$\epsilon_q'$ be as in Example 3.7. Then any CM abelian variety
over $K$ of type $\epsilon'=\epsilon \epsilon_q'$ is of dimension
$2h$.
 \endexample

\subheading{4. Descent to $\Bbb Q$}
%Canonical Hecke characters and proof of theorem 2}

 Let $\chi$ be a Hecke character of $K$ satisfying $(3.1)$ and let
$A(\chi)$ be an associated abelian variety over $K$ of CM type
$(T, \Phi)$ with $T=\Bbb Q(\chi)$ and $\Phi$ being the set of
complex embeddings of $T$ which are the identity on $K$. If
$A(\chi)$ descends to an abelian variety over $\Bbb Q$, then
$(0.2)$ holds. For convenience, we repeat it here as
$$
\chi(\bar\frak a) =\overline{\chi(\frak a)} \tag{4.1}
$$
for every ideal of $K$ prime to the conductor $\frak f$ of $\chi$.
Conversely,  the proof of [Sh1, Proposition 5, Page 521] gives

\proclaim{Lemma 4.1} Let $\chi$ be a Hecke character of $K$
satisfying $(3.1)$ and $(4.1)$. Then there is a CM abelian variety
$(A, i)$ over $K$ of CM type $(T, \Phi)$ associated to $\chi$ such
that $(A, i_{T^+})$ is actually defined over $\Bbb Q$. Here
$T=\Bbb Q(\chi)$, $\Phi$ is the set of embeddings of $T$ which is
the identity on $K$, $T^+$ is the maximal totally real subfield of
$T$ and $i_{T^+}$ is the restriction of $i$ on $T^+$.
\endproclaim
\demo{Proof} Choose $\delta \in T^*$ such that
$\bar\delta=-\delta$. Then there is, by [Sh1, Theorem 6, Page
512], a structure $(A, \Cal C, i)$ of type $(T, \Phi, \O_T,
\delta)$ rational over $K$ which determines $\chi$. Here $\Cal C$
is a polarization of $A$. We refer to [Sh1, Page 509] for the
meaning of type $(T, \Phi, \O_T, \delta)$.  Let $\rho$ be the
complex conjugation of $\Bbb C$, restricting to $T$ or $K$. Then
[Sh1, Lemma 3, Page 520] asserts that $(A^\rho, \Cal C^\rho, i^*)$
is of the  same type $(T, \Phi, \O_T, \delta)$, where $i^*(a)
=i(a^\rho)^\rho$. So there is an isomorphism $\mu$ from $(A, \Cal
C, i)$  to  $(A^\rho, \Cal C^\rho, i^*)$ over $\Bbb C$. On the
other hand, $(4.1)$ implies that they determine the same Hecke
character $\chi$ by [Sh1, Proposition 1, Page 511], and thus any
isogeny between the two structures is defined over $K$. In
particular, the isomorphism $\mu$ is defined over $K$. Let
$$
\omega: T_{\Bbb R} = T\otimes_{\Bbb Q} \Bbb R \lr A
$$
be an homomorphism defining $(A, \Bbb C, i)$ as of type $(T, \Phi,
\O_T, \delta)$, and let $\omega'(u)= \omega(\bar u)^\rho$. Then
$\omega'$ is an homomorphism from $T_{\Bbb R}$ to $A^\rho$ by
[Sh1, Lemma 3, Page 520], and $\mu$ can be chosen so that
$\mu'=\mu\circ \omega$. So
$$
\omega(u) =\omega'(\bar u)^\rho = \mu^\rho(\omega(\bar)^\rho)
   =\mu^\rho \mu (\omega(u))
   $$
   and thus $\mu^\rho \mu=1$. By the descent theory, $(A, \Cal C, i_{T^+})$ can thus
   be descended to $\Bbb Q$.
   \enddemo

By [Ro2, Proposition 1], the condition $(4.1)$ is equivalent to
each of the two following conditions
$$
\epsilon(n) = \left(\frac{-D}n \right) \tag{4.2}
$$
when $n$ is prime to $\frak f$, or
$$
\chi^{un}|_{\Bbb Q_{\A}^*} =\kappa  .\tag{4.3} $$ Here
$\kappa=\prod \kappa_p$ is the quadratic character of the ideles
$\Bbb Q_\A^*$ associated to $K/\Bbb Q$ by the global class field
theory, and $\chi^{un} =\chi |\,|_{\A}^{\frac{1}2}$ is the
unitarization of $\chi$,  viewed as an idele class character of
$K$.
 In
particular, $\frak f$ is divisible by every prime ideal of $K$
which is ramified in $K/\Bbb Q$. In this section, we consider the
set $E$ of characters
$$
\epsilon: \left(\O_K/\frak f\right)^* \lr \mu_{12} \tag{4.4}
$$
satisfying  $(4.2)$ and that $\frak f$ is only divisible by
ramified primes of $K$. For simplicity, we assume $D >4$ so
$\O_K^* =\{ \pm 1 \}$. For $\epsilon \in E$, we will determine the
dimension and conductor of a CM abelian variety over $K$ of type
$\epsilon$. Write
 $$
 \frak f = \prod_{p|D} \frak p^{e_p}
 $$
 then
 $$
 \left(\O_K/\frak f\right)^* \cong \oplus_{p|D} (\O_\frak p/\frak
 p^{e_p})^*
.
$$
Here $\frak p$ is the unique prime ideal of $K$ above $p$. So a
character $\epsilon$ satisfying $(4.4)$ is the same as a character
$$
\epsilon=\prod_{p|D} \epsilon_p: \, \prod_{p|D} \O_p^* \lr
 \mu_{12} \tag{4.4'}
$$
with $\epsilon_p: \O_p^* \lr \mu_{12}$. Here $K_p=K\otimes \Bbb
Q_p$ is the completion of $K$ at the prime above $p$, and $\O_p$
is the ring of integers in $K_p$.  As remarked in [Ro2, Page 522],
$\chi_p^{un} = \epsilon_p^{-1}$ on $\O_p^*$, a fact we will use
later. In particular, that $\epsilon$ satisfies $(4.2)$ is
equivalent to that every $\epsilon_p$ satisfies
$$
\epsilon_p|_{\Bbb Z_p^*}= \kappa_p. \tag{4.5}
$$
When $ p \nmid 6$, $(\O_p)^* \cong (\O_p/\frak p)^* \times (1+
\varpi_p \O_p)$, and $(\O_p/\frak p)^* \cong (\Bbb Z/p)^*$, where
$\varpi_p$ is uniformizer of $\O_p$.  So there is a unique
character $\epsilon_p:\O_p^* \lr \mu_{12}$ satisfying $(4.5)$, and
it is trivial on $1+\varpi_p \O_p$. We denote it by
$\epsilon_p^0$.  Similarly, when $p=3$, there is a unique
character $\epsilon_3^0$ of $\O_3^*$ of order $2$ and conductor
index $1$ satisfying $(4.5)$. Here the conductor index of a
character $\chi$ of $K_p^*$ (or $\O_p^*$) is the smallest integer
$r \ge 0$ such that $\chi|_{1+\varpi_p^r\O_p}=1$. Since
$1+\varpi_3\O_3$ is  a cyclic pro-3-group, there are also two
characters of $\O_3^*$ of order $6$ and conductor index $2$
satisfying $(4.5)$, given by $\epsilon_3^0 \phi_3^{\pm 1}$, where
$\phi_3$ is trivial on $(\O_3/\frak p_3)^*= \{\pm 1\}$ and of
order $3$ on $1+\varpi_3 \O_3$. So the set $E_3$ of characters of
$\O_3^*$ of order $\le 12$ satisfying $(4.5)$ is
$$
E_3=\{\epsilon_3^0 \phi_3^i:\, -1 \le i \le 1\}. \tag{4.6}
$$

 The case $p=2$ is a little more complicated and interesting.
Let $G=\{ z \in K_2: z \bar z =1\}$ be the norm one group, let
$$
G_n=\{ z= x + y \sqrt{-d}:\, y \in 2^n \Bbb Z_2, x-1 \in 2^n d
\Bbb Z_2\}
$$ be subgroups of $G$ for $n \ge 0$. Here $d =D/4$. Then $G=G_0$
and $[G_n: G_{n+1}]=2$. One has an exact sequence
$$
1 \lr \Bbb Z_2^* \lr \O_2^* \lr G_1 \lr 1.
$$
The last map is $z \mapsto  z/\bar z$.

When $8|D$, $2|d$,  It is known ([Ro2, Proposition 4]) that
$\kappa_2$ extends to two quadratic characters $\epsilon_2^{\pm}$
of $\O_2^*$ of conductor index $5$. Since $G_1$  is a cyclic
pro-2-group, the set $E_2$ of characters of $\O_2^*$ of order $\le
12$ satisfying $(4.5)$ is
$$
E_2 =\{ \epsilon_2^\pm , \epsilon_2^\pm
\tilde\phi\}=\{\epsilon_2^+ \tilde\phi^i:, -1 \le i \le 2\}.
\tag{4.7}
$$
Here $\phi$ is a fixed character of $G_1$ of order $4$, trivial on
$G_3$, and $\tilde\phi(z) =\phi(z/\bar z)$.  Since
$$
1 \lr \Bbb (Z/8)^* \lr (\O_2/\varpi_2^5)^* \lr G_1/G_3 \lr 1,
$$
all characters in $E_2$ have conductor index $5$.

  When $4||D$, $\sqrt{-d}$ is a unit in $\O_2$ and
  $\varpi_2=1+\sqrt{-d}$ is a uniformizer of $K_2$. In this case,
  $G_1= \{ \pm 1 \} \times G_2$ and $G_2$ is a cyclic pro-2-group.
  It is easy to check that $(\O_2/\varpi^3)^*$ is cyclic of order
  $4$, generated by $\delta= \sqrt{-d}$.  So
  $$
  \epsilon_2^\pm: \O_2^* \lr (\O_2/\varpi_2^3)^* \lr \mu_4,
   \sqrt{-d} \mapsto \pm i
   \tag{4.8}
   $$
   gives two extensions of $\kappa_2$ to $\O_2^*$. Let $\phi$ be a
   fixed character of $G_1$ of order $4$, trivial on $\{\pm 1\} \times G_4$. Then
   the set $E_2$ of characters of $\O_2^*$ of order $\le  12$
   satisfying $(4.5)$ is
   $$
   E_2=\{ \epsilon^\pm \tilde\phi^i: -1 \le i \le 2\}. \tag{4.9}
   $$
  The characters $\epsilon_2^\pm \tilde\phi^i$ have conductor
  indices $3$, $5$, and $7$ respectively according to $i=0, 2, \pm
  1$. So we have  the following proposition.

  \proclaim{Proposition 4.2} Assume $D > 4$. Let $E$ be the set of characters of
  $(O_K/\frak f)^*$ of order $\le 12$ satisfying $(4.2)$ and that
  $\frak f$ is only divisible by primes dividing $\sqrt{-D}\O_K$.
  Then
  $
  \# E = \# E_2 \cdot\# E_3,
  $
  where $E_p$ is the set defined in $(4.6)$, $(4.7)$, or $(4.9)$
  when $p |D$ and trivial otherwise. Moreover

  (1) \quad When $3|D$, $\# E_3 = 3$, one of which has conductor
  index $1$, and the other two have conductor index $2$.

   (2) \quad When $8|D$, $\# E_2 = 4$, all have conductor index
   $5$.

    (3) \quad When $ 4||D$, $\# E_2 =8$, and their
    conductor indices are $3$, $5$, or $7$.

    (4) \quad Let $\epsilon =\prod_{p|D} \epsilon_p \in E$ and let $\chi$
    be a Hecke character of $K$ of type  $\epsilon$.
    The the conductor of $\chi$ is $\prod_{p|D} \frak p^{e_p}$
    where $e_p$ is the conductor index of $\epsilon_p$. One has
    $e_p=1$ for $p\ne 6$.

    (5) \quad Let  $\epsilon \in E$, let $\chi$
    be a Hecke character of $K$ of type  $\epsilon$,   and let $A$ be a CM abelian
    variety over $K$ associated with $\chi$. Then $\dim A= h$ unless
    $4|D$ and $\epsilon_2 =\epsilon_2^\pm \tilde\phi$, where $\phi$
     is a
    fixed character of $G_1$ of order $4$ given above.  In such a
    case,
    $$
    \dim A = \cases
         h &\ff d=D/4 \equiv 1 \mod 8,
         \\
          2h &\hbox{ otherwise}.
          \endcases
          $$
          \endproclaim
  \demo{Proof} Every claim except (5) was proved above. We use
  Theorem 3.5 to verify (5) in the case $2|D$, and leave the other
  cases to the reader.

    We first assume $8|D$ so that $2|| d$. If
  $\epsilon$ is quadratic, it is obvious by Theorem 3.5. If
  $\epsilon_2 =\epsilon_2^\pm$ is quadratic, and
  $\epsilon_3 =\epsilon_3^0 \phi_3 \in E_3$ is of order $6$ as in $(4.6)$.
   Let $\alpha_3 = \frac{\beta}{\sqrt{-D}}
  \in K^*$ be such that $3 \alpha_3^2$ is prime to $D$ and that
  $\beta$ is prime to $3$. Then
  $$
  \align
  \epsilon(3 \alpha_3^2)
    &= \phi_3(3 \alpha_3^2) \epsilon_2^\pm (3)  \epsilon_3^0(-3/D)
       \prod_{p|D, p\nmid 6} \epsilon_p^0(3)
       \\
       &= \phi_3(3 \alpha_3^2) (-3/D, -D)_3 \prod_{p \ne 3} (3, -D)_p
       \\
        &=\phi_3(3 \alpha_3^2) \prod_{p<\infty} (3, -D)_p (-1,
        -D)_3
        \\
        &=-\phi_3(3 \alpha_3^2)
        \endalign
        $$
        generates $\hbox{Im}(\epsilon)$.
        So Theorem 3.5(3) implies that $\dim A =h$.
        When $\epsilon_2 =\epsilon_2^\pm \tilde\phi$ is of order $4$,
        let $\alpha_2 = \frac{\beta}{\sqrt{-d}}$ such that $2
        \alpha_2^2$ is prime to $D$ and  $\beta$ is prime to
        $2$. Then
        $$
        \align
         \epsilon(2 \alpha_2^2) &=\pm  \epsilon_2(2 \alpha_2^2)\prod_{2\ne
         p|D} \epsilon_p(2 \alpha_2^2)
         \\
          &=\pm  \kappa_2(-2/d) \prod_{2 \ne p |D} \epsilon_p(2 \alpha_2^2)
          \endalign
         $$
         is of order $2$ or $6$. So Theorem $3.5$ asserts
         that $\dim A =2h$.

          Finally, assume $4||D$, let $\alpha
         =\frac{1+\sqrt{-d}}2 \beta \in K^*$ be such that
         $2 \alpha_2^2 = (\frac{1-d}2 +\sqrt{-d}) \beta^2$ is
         prime to $D$ and $\beta$ is prime to $2$. For
         $\epsilon_2= \epsilon_2^\pm  \tilde\phi^j$ with $ -1 \le j \le
         2$.  Then
         $$
 \align
         \epsilon(2 \alpha_2^2)
          &=\pm  \epsilon_2^{\pm}(\frac{1-d}2 +\sqrt{-d})
               \phi^j (g) \prod_{2 \ne p|D}\epsilon_p(2 \alpha_2^2)
               \\
               &=\pm i \phi^j (g) \prod_{2 \ne p|D}\epsilon_p(2 \alpha_2^2)
               \endalign
               $$
               is of order $4$ or $12$ if and
               only if $\phi^j(g) =\pm 1$. Here
               $$
                g=\frac{\sqrt{-d} + \frac{1-d}2}{\sqrt{-d}
                -\frac{1-d}2} \in G_2.
                $$
 When $j=0$ or $2$, one  has always $\phi^j(g) =\pm 1$. When $
 j=\pm 1$, $\phi^j(g) =\pm 1$ if and only if $g \in G_3$, which is
 in turn equivalent to $d \equiv 1 \mod 8$. This finishes the
 proof.
 \enddemo

\proclaim{Corollary 4.3} Assume $D >4$.  Let the notation be as in
Proposition 4.2. Let $E_{Sim}$ be the subset of characters
$\epsilon$ in $E$ such that a CM abelian variety $A$  over $K$ of
type $\epsilon$ has dimension $h$ and the smallest conductor. Then
$$
E_{Sim} =\cases
     \{ \prod_{p|D} \epsilon_p^0 \} &\ff  2 \nmid D,
     \\
      \{ \epsilon_2^\pm \prod_{2\ne p|D} \epsilon_p^0 \} &\ff 2|D.
      \endcases
      $$
      In particular, $\#E_{Sim}=1$ or $2$ depending on whether $D$
      is odd or even.
      \endproclaim

  {\bf Proof of Theorem 0.2} We may assume $D>4$. By Remark 1.2,  a simple CM
  abelian variety $A$ over $K$  is up to isogeny one-to-one
  correspondence to the set $\{ \chi^{\sigma}: \sigma \in
  \hbox{Aut}(\Bbb C/K)\}$, where $\chi$ is a Hecke character of
  $K$ associated to $(A, i)$ of CM type $(T, \Phi)$ for some embedding
  $i$ and some type $(T, \Phi)$. Here we require all the associated
  Hecke characters of $K$ to be of infinite type $\{\hbox{ id }\}$.
  When $D$ is odd or divisible by $8$, $\epsilon$ is quadratic,
  and
  $$
  \{\chi^\sigma: \, \sigma \in \hbox{Aut}(\Bbb C/K)\}
   = \{ \chi \phi: \, \phi \hbox{ is an ideal class character of } K \}
   $$
   is uniquely determined by $\epsilon$. So Theorem 0.2 follows
   from Corollary 4.3.  Assume now that
   $4||D$ and $D>4$, and  let $L=K(i)$.  Fix  $\rho' \in
   \hbox{Aut}(\Bbb C/K)$ such that $\rho'(i) =-i$, and an ideal
   class character $\phi'$ of $K$ such that $\phi'(\frak p_2) =-1$
   for the prime ideal $\frak p_2$ above $2$.  Fix an character $\epsilon$
   in Corollary $4.3$ and a Hecke character $\chi_0$ of $K$
    of type $\epsilon$.  Then  the other character in Corollary 4.3 is
   $\epsilon^{-1}$, and $\chi_0'=\chi_0^{\rho'}$ is of type $\epsilon^{-1}$.
    The Hecke characters of $K$ satisfying $(0.2)$ such that its
    associated abelian varieties have the smallest conductor and
    dimension  are of type $\epsilon^{\pm 1}$ by Corollary 4.3 and
    are given by
    $$
    \align
    &\{ \chi_0 \phi, \chi_0' \phi:\,
    \phi\hbox{ is an ideal class character of } K \}
    \\
     &=\{\chi_0 \phi, \chi_0' \phi:\, \phi\hbox{ is an ideal class character of   $K$
     such that } \phi(\frak p_2) =1 \}
     \\
      &\quad \bigcup
          \{\chi_0 \phi' \phi, \chi_0'\phi' \phi:\, \phi\hbox{ is an ideal class character of   $K$
     such that } \phi(\frak p_2) =1 \}
     \\
      &= \{ \chi_0^\sigma: \sigma \in \hbox{Aut}(\Bbb C/K) \}
         \bigcup  \{ (\chi_0\phi')^\sigma: \sigma \in \hbox{Aut}(\Bbb C/K)
         \}.
         \endalign
         $$
         Here we used the fact that
         $$
         \align
         &\{ \chi \phi:\, \phi\hbox{ is an ideal class character of   $K$
     such that } \phi(\frak p_2) =1 \}
     \\
       &= \{ \chi^\sigma:\, \sigma \in \hbox{Aut}(\Bbb C/L)\}.
       \endalign
       $$
         Therefore, there are exactly two CM abelian varieties
         over $K$ of dimension $h$, up to isogeny, which descend
         to $\Bbb Q$ and have the smallest conductor.

  From the proof, it is clear that the isogeny class of
  the CM abelian varieties $A$ are not determined by $\epsilon$ when
  the image of $\epsilon$ is not in $K$.

\subheading{5. Root number}

  Let $\chi$ be a simplest Hecke character of $K$ of type
  $\epsilon$, given by Corollary 4.3. The condition $(0.2)$
  implies $W(\chi)= \pm 1$. It is known ([Ro3]) that
  $$
  W(\chi) =\cases
     \left(\frac{2}{D}\right) &\ff  2 \nmid D,
     \\
      \epsilon(1+\sqrt{-D/4}) &\ff 8||D.
      \endcases
      \tag{5.0}
      $$
 We may thus assume that $4||D$ and $D> 4$.  Let $\frak p_2$ be the prime ideal of $K$
  above $2$, $d=D/4$, and $\alpha_0=1+\sqrt{-d}$. Then $\epsilon$
  can be viewed as
  $$
   \epsilon=\epsilon_2 \epsilon^0: \left(\O_K/\frak p_2^3
   \sqrt{-d}\right)^* = \left(\O_K/\frak p_2^3\right)^* \times
   \left(\O_K/\sqrt{-d}\right)^* \lr \mu_4
   \tag{5.1}
$$
 Here
$\epsilon_2$ is trivial on $(\O_K/\sqrt{-d})^*$ but gives an
isomorphism $\left(\O_K/\frak p_1^3\right)^* \cong \mu_4$, while
$\epsilon^0$ is trivial on $\left(\O_K/\frak p_2^3\right)^*$ and
is $(\frac{}d)$ on $\left(\O_K/\sqrt{-d}\right)^* \cong (\Bbb
Z/d)^*$.

\proclaim{Proposition 5.1} Let the notation be as above, and let
$\chi $ be a simplest Hecke character of $K=\Bbb Q(\sqrt{-D})$ of
type $\epsilon$.  Then
$$
\chi(\alpha_0^{-1}\frak p_2) =W(\chi) (1-\epsilon_2(\sqrt{-d}))
\alpha_0^{-1}.
$$
\endproclaim
\demo{Proof}
 For the proof, we follow [Ro2] closely. Let $\kappa=\prod_p
 \kappa_p$ and  $\chi^{un} =\chi |\,|_{\A}^{\frac{1}2}=\prod_p \chi_p^{un}$ be as in
 as in Section 4. Then
 $$
 \chi^{un}|_{\Bbb Q_\A^*} =\kappa,  \quad \hbox{ and }
  \quad \chi_p^{un}|_{\O_p^*} =\epsilon_p^{-1}  \hbox{ for }
  p|D.
  \tag{5.1}ite
  $$
  Here we write $\epsilon=\prod_{p|D}\epsilon_p$  as in Corollary 4.3.
  Let $\psi=\prod\psi_p$ be a nontrivial additive character of
  $\Bbb Q_\A$ and let $\psi_K=\psi \circ N_{K/\Bbb Q}$. Then
  Rohrlich proved in [Ro2] that the relative local root number
  $$
  W(\kappa_p, \chi_p) = \frac{W(\chi_p, \psi_{K_p})}{W(\kappa_p,
  \psi_p)}
  $$
  is independent of the choice of $\psi$. Here $W(\chi_p,
  \psi_{K_p})$ and $W(\kappa_p, \psi_p)$ are local root numbers.
  Since the global root number of $\kappa$ is one, we have then
  $$
  W(\chi) =\prod_{p\le \infty} W(\kappa_p, \chi_p).
  \tag{5.2}
  $$
  By [Ro2, Propositions 8, 11, and 12], one has
  $$
  W(\kappa_p, \chi_p) =\cases
     1 &\ff p\nmid D,
     \\
      \kappa_p(2) &\ff 2\ne p|D.
      \endcases
      \tag{5.3}
      $$
So
$$
W(\chi) =\prod_{2\ne p|D}  \kappa_p(2) W(\kappa_2,
\chi_2)=\kappa_2(2) W(\kappa_2, \chi_2). \tag{5.4}
$$
On the other hand, since $\alpha_0^{-1}\frak p_2$ is prime to
$D\O_K$, one has
$$
\align
 \chi(\alpha_0^{-1}\frak p_2)&=\prod_{p\nmid D\infty}
 \chi_p(\alpha_0^{-1})
 \\
 &=\prod_{p|D\infty} \chi_p(\alpha_0)
 \\
 &=\chi_2(\alpha_0) \alpha_0^{-1}.
 \endalign
 $$
 Combining this with $(5.4)$ and the following lemma, one proves
 the proposition.
\enddemo

\proclaim{Lemma 5.2} Let the notation be as above. Then
$$
W(\kappa_2, \chi_2) =
\frac{1+\chi_2(\sqrt{-d})}{\chi_2(1+\sqrt{-d})}\kappa_2(2) .
$$
\endproclaim
\demo{Proof} We first remark that $\chi_2(\sqrt{-d}) =
\epsilon_2(\sqrt{-d})^{-1}= - \epsilon_2(\sqrt{-d})$ by $(5.1)$.
  We
recall a general way to compute local root numbers. Let $F$ be a
non-archimedean local field with ring of integers $\O_F$ and a
uniformizer $\varpi$. Let $\theta$ be a  unitary character of
$F^*$ of conductor index $r$ and let $\psi_F$ be a nontrivial
additive character of $F$ of conductor index $n=n(\psi_F)$, i.e.,
the smallest integer $n$ such that $\psi|_{\varpi^n\O_F}=1$. Then
$$
W(\theta, \phi_F) = b^{-1} \theta^{-1}(\beta)
    \int_{\O_F^*} \theta^{-1}(x) \psi_F(x \beta) dx
    \tag{5.5}
    $$
    where  $\beta \in \varpi^{-r+n} \O_F^*$ and
    $b = |\varpi|^{\frac{r}2} \hbox{Vol}(\O_F, dx)$.
Now choose $\psi_2 : \Bbb Q_2 \lr \Bbb C^*$ via  $\psi_2(x) = e^{2
\pi i \lambda_2(x)}$ where $\lambda_2$ is the canonical map
$$
\lambda_2: \Bbb Q_2 \lr \Bbb Q_2/\Bbb Z_2 \hookrightarrow \Bbb
Q/\Bbb Z.
$$
Then a simple calculation gives
$$
W(\kappa_2, \psi_2) = \cases
    i &\ff 4||D,
    \\
     i \kappa_2(2) &\ff 8||D \hbox{ and } \kappa_2(-1)=-1,
     \\
      \kappa_2(2) &\ff  8||D \hbox{ and } \kappa_2(-1)=1.
     \endcases
     \tag{5.6}
     $$
As for $W(\chi_2, \psi_{K_2})$, one has $\chi_2^{un}(a) =
\chi_2(a) |a|_2^{\frac{1}2}$, and in particular,
$\chi_2^{un}(\alpha_0) = \chi_2(\alpha_0)/ \sqrt{2}$.
 Since the conductor index of $\chi_2$ and
$\psi_{K_2}$ are $3$ and $-2$ respectively, and $\alpha_0$ is a
uniformizer of $K_2$,  we may take $\beta =\alpha_0/8$ in $(5.5)$
and have
$$
\align
 W(\chi_2, \psi_{K_2})&= W(\chi_2^{un}, \psi_{K_2})
  \\
   &= \frac{2\sqrt 2}{\hbox{Vol}(O_2, dz)} \chi_2^{un}
   (\alpha_0/8)^{-1}
     \int_{O_2^*} \chi_2^{un}(z^{-1}) \psi_{K_2}(z \alpha_0/8) dz
     \\
     &= \frac{1}{2 \sqrt 2}\kappa_2(2) \chi_2^{un}(\alpha_0)^{-1}
       \sum_{a\in (O_2/\frak p_2^3)^*} \chi_2^{un}(a)^{-1} \psi_2(
       \frac{\tr a \alpha_0 }{8}).
       \endalign
       $$
       Since $(O_2/\frak p_2^3)^* =\{ \pm 1, \pm \sqrt{-d}\}$, one
       has
       $$
       \align
    W(\chi_2, \psi_{K_2})
    &=i \kappa_2(2) \chi_2^{un}(\alpha_0)^{-1}
    \frac{1+\chi_2^{un}(\sqrt{-d})}{\sqrt{2}}
    \\
      &=W(\kappa_2, \psi_2) \kappa_2(2) \frac{1+\chi_2(\sqrt{-d})}{\chi_2(1+\sqrt{-d})}.
    \endalign
    $$
\enddemo

\subheading{6. The central $L$-value}

  The purpose of this section is to prove the following theorem.

  \proclaim{Theorem 6.1} Let $\chi$ be a simplest Hecke character
  of $K=\Bbb Q(\sqrt{-D})$. Then the central $L$-value
  $L(1, \chi)\ne  0$ if and only if $W(\chi)= 1$.
  \endproclaim

   The case where $D$ is odd or divisible by $8$ was proved by
   Montgomery and Rohrlich ([MR]) in 1982. We will use the same
   method to settle the case $4||D$ and $D>4$. The case $D=4$ is well-known.
    We assume again in this section
   $4||D$.  One complication is that
   the two family $\{\chi^\sigma:\, \sigma \in
   C/K)\}$ and $\{ \chi \phi:\,  \phi  \hbox{ is an ideal class
   character of } K\}$ are not the same in our case. Let $L=K(i)$ and fix
   $\rho' \in \hbox{Aut}(\Bbb C/K)\}$ such that $\rho'(i) = -i$.
   If $\chi$ is of type $\epsilon$ where $\epsilon=\epsilon_2
   \epsilon^0$ is given by $(5.1)$, then $\chi'=\chi^{\rho'}$ is
   of type $\epsilon^{-1}= \epsilon_2^{-1} \epsilon^0$. For each
   ideal class $C$ of $K$, let
   $$
   L(s, C, \chi) =\sum_{\frak a \in  \Cal C, \hbox{ integral} }
   \chi(\frak a) (N \frak a)^{-s}
   \tag{6.1}
   $$
   be the partial L-function, and denote
   $$\tilde L(s, C, \chi)= L(s, C, \chi) + L(s, C, \chi').
   \tag{6.2}
   $$
   We also write $[\frak a]$ for the ideal class of $K$ containing
   the ideal $\frak a$. Since
   $$
   \{ \chi^\sigma:\, \sigma \in \hbox{Aut}(\Bbb C/L)\}
   = \{ \chi \phi:\, \phi  \hbox{ is an ideal class
   character of  $K$ with } \phi([\frak p_2]) =1\},
   $$
   one has by the same argument as in [Ro3, Page 226]

   \proclaim{Lemma 6.2} Let the notation be as above. Then the
   following are equivalent.

   $(1)$ \quad The central L-value $L(1, \chi) = 0$,

    $(2)$ \quad The partial central  L-values
     $L(1, C, \chi) + L(1, C[\frak p_2], \chi)=0$ for every ideal
     class $C$ of $K$.

    $(3)$ \quad The partial central  L-values
     $\tilde L(1, C, \chi) + \tilde L(1,  C[\frak p_2], \chi)=0$ for every ideal
     class $C$ of $K$.
     \endproclaim

 Following [MR], we recall that the theta function
 $$
 \theta(t) = h + 2 \sum_{n \ge 1} a(n) e^{-\frac{2 \pi n}{\sqrt
 D}}\tag{6.3}
 $$
 is strictly increasing  for $ t >0$ and satisfies the simple functional
 equation
 $$
 \theta(1/t) =t \theta(t). \tag{6.4}
 $$
 Here $a(n)$ is the number of integral ideals of $K$ with norm
 $n$. Recall also that the Eisenstein series
 $$
 G(z, s) = \sum_{n >0, m \in \Bbb Z} \left(\frac{-D}n\right) (m D
 \bar z +n) |m Dz + n|^{-2s}\tag{6.5}
 $$
 has analytic continuation for all $s \in \Bbb C$ and $z$ in the
 upper half plane, and
 $$
  \frac{\sqrt D}{\pi} G(z, 1)
   =h + 2 \sum_{n \ge 1} a(n) e^{2 \pi n i z}.
  $$
In particular,
$$
\frac{\sqrt D}{\pi } G(\frac{i t}{\sqrt D}, 1) = \theta(t).
\tag{6.6}
$$
It is convenient to denote
$$
G_{\odd}(z, s)=G(z, s) -G(2z, s)
    =\sum_{m \, \odd, n >0} \left(\frac{-D}n\right) (m D
 \bar z +n) |m Dz + n|^{-2s}.
 $$
In particular,
$$
\frac{\sqrt D}{\pi} G_{\odd}(\frac{i t}{\sqrt D}, 1) =
\theta(t)-\theta(2t) >0. \tag{6.7}
$$

Since we can switch between $\chi$ and $\chi'$, we may and will
assume that $\epsilon_2(\sqrt{-d})=i$.

\proclaim{Lemma 6.2} Let $\epsilon=\epsilon_2 \epsilon^0$ with
$\epsilon_2(\sqrt{-d})=i$.

(1) \quad When  $\beta = m \sqrt{-d} +n$  is prime to $D\O_K$, $m,
n \in \Bbb Z$,  one has
$$
\epsilon(\beta) =\cases
  \left(\frac{-D}n\right) &\ff m \equiv 0 \mod 4,
  \\
 -\left(\frac{-D}n\right) &\ff m \equiv 2 \mod 4,
 \\
  i \left(\frac{-1}m\right) \left(\frac{n}d\right) &\ff n \equiv 0
  \mod 4,
  \\
   - i \left(\frac{-1}m\right) \left(\frac{n}d\right) &\ff n \equiv
   2
  \mod 4.
 \endcases
 $$

 (2) \quad When $\beta =\alpha_0 (a + b \frac{1-\sqrt{-d}}2) \in
 \alpha_0 \frak p_2^{-1}$  is prime to $D \O_K$, where $\alpha_0
 =1+\sqrt{-d}$, $a, b \in \Bbb  Z$, one has
 $$
 \epsilon(\beta) = \cases
   \left(\frac{-D}{2a+b}\right) &\ff a \equiv 0 \mod 2,
   \\
    i \left(\frac{-D}{2a+b}\right) &\ff a \equiv 1 \mod 2.
    \endcases
    $$
    \endproclaim
\demo{Proof} $(1)$ follows from definition easily. For example, if
$m \equiv 2 \mod 4$, then $\beta =n-m + m\alpha_0 \equiv  n-m \mod
\frak p_2^3$ and thus
$$
\epsilon(\beta)
 = \epsilon_2(\beta) \left(\frac{n}d\right)
= \left(\frac{-1}{n-m}\right) \left(\frac{d}n\right)
 =-\left(\frac{-D}n\right).
   $$
 (2) \quad Notice that
   $$
   \beta = a+ \frac{1+d}2 b +a \sqrt{-d}.
   $$
  When $a \equiv 0 \mod 4$, one has
   $$
   \align
   \epsilon(\beta) &= \left(\frac{-1}{a+\frac{1+d}2 b}\right)
                       \left(\frac{a + \frac{1+d}2 b} d\right)
   \\
    &=\left(\frac{2}d\right) \left(\frac{2a+b}d\right)
       \left(\frac{-1}{a+\frac{1+d}2 b}\right)
       \\
       &=\cases
            \left(\frac{d}{2a+b}\right) \left(\frac{-1}b\right)
                       &\ff d \equiv 1 \mod 8
    \\
        - \left(\frac{d}{2a+b}\right) \left(\frac{-1}{3b}\right)
                       &\ff d \equiv 5 \mod 8
                        \endcases
   \\
    &=\left(\frac{-D}{2a +b}\right).
    \endalign
    $$
    The case $a \equiv 2 \mod 4$ is similar. When $a$ is odd, and
    $d \equiv 1 \mod 8$, \newline
    \noindent $a +\frac{1+d}2 b \equiv a+b \mod 4$, and
    so
    $$
    \align
     \epsilon(\beta) &=
        \cases i \left(\frac{-1}a\right) \left(\frac{2a+b}d
        \right) \left(\frac{2}d\right) &\ff a+b \equiv 0 \mod 4,
        \\
         - i \left(\frac{-1}a\right) \left(\frac{2a+b}d
        \right) \left(\frac{2}d\right) &\ff a+b \equiv 2 \mod 4
         \endcases
    \\
     &=i \left(\frac{-1}{2a+b}\right) \left(\frac{d}{2a+b}\right)
     \\
      &= i \left(\frac{-D}{2a+b}\right).
      \endalign
      $$
      The case when $a$ is odd and $d \equiv 5 \mod 8$ is similar.
\enddemo

\proclaim{Lemma 6.3} One has
$$
\tilde L(s, [\O_K], \chi) = 4 G(\frac{2i}{\sqrt D}, s)-2
G(\frac{i}{\sqrt D}, s), \tag{6.8}
$$
and
$$
\tilde L(s, [\frak p_2], \chi)
 = 2^s W(\chi) G_{\odd} (\frac{i}{2\sqrt D}, s ).
 \tag{6.9}
 $$
 \endproclaim
 \demo{Proof} Since $D>4$, $\O_K^* =\{\pm 1\}$, one has by Lemma
 6.2(1)
 $$
 \align
 \tilde L(s, [\O_K], \chi)
 &= \sum_{\alpha \O_K} (\chi(\alpha \O_K) + \chi'(\alpha\O_K))
 (N\alpha)^{-s}
 \\
  &= \sum_{\alpha \in \O_K} \frac{\epsilon(\alpha) +
  \epsilon^{-1}(\alpha)}2 \alpha (N\alpha)^{-s}
  \\
  &= 2 \sum_{m\equiv 0(4), n>0} \left(\frac{-D}n\right) (m
  \sqrt{-d} +n) |m \sqrt{-d} +n|^{-2s}
  \\
   &\qquad
   -2 \sum_{m\equiv 2(4), n>0} \left(\frac{-D}n\right) (m
  \sqrt{-d} +n) |m \sqrt{-d} +n|^{-2s}
  \\
   &=2 G(\frac{2i}{\sqrt D}, s)
     - 2 G_{\odd} (\frac{i}{\sqrt D}, s)
     \\
      &= 4 G(\frac{2i}{\sqrt D}, s)- 2 G (\frac{i}{\sqrt D}, s).
      \endalign
      $$
      This proves $(6.8)$.  On the other hand, $\frak a \in [\frak
      p_2]$ is integral and prime to $D\O_K$ if and only if
      $ \frak a= \beta \alpha_0^{-1} \frak p_2$ with $\beta \in
      \alpha_0 \frak p_2^{-1}$ prime to $D \O_K$. In such a case,
      Proposition 5.1 gives
      $$
      \align
      \chi(\frak a)
       &=\epsilon(\beta) \beta \chi(\alpha_0^{-1} \frak p_2)
       \\
        &=W(\chi) (1-\epsilon_2(\sqrt{-d})) \alpha_0^{-1}
         \epsilon(\beta) \beta.
       \endalign
        $$
     One also has  $W(\chi) =W(\chi')$ by Proposition 5.1. So
        $$
        \align
         \chi(\frak a) + \chi'(\frak a)
         &=W(\chi) \alpha_0^{-1} \beta \left[
          \epsilon(\beta) (1 -\epsilon_2(\sqrt{-d}))
           +\epsilon(\beta)^{-1} (1 + \epsilon_2(\sqrt{-d}))\right]
           \\
           &=W(\chi) 2 \alpha_0^{-1} \beta \cases
              \epsilon(\beta) &\ff \epsilon(\beta) =\pm 1,
              \\
               - \epsilon(\beta) \epsilon_2(\sqrt{-d})
                     &\ff \epsilon(\beta) =\pm i.
                     \endcases
        \\
         &=W(\chi)(2a+b-b\sqrt{-d}) \left(\frac{-D}{2a+b}\right),
          \endalign
          $$
 if  $\beta =\alpha_0 (a +b \frac{1 -\sqrt{-d}}2)
 \in \alpha_2 \frak p_2^{-1}$  as in Lemma 6.2, and
 $\epsilon_2(\sqrt{-d})=i$ (our assumption). So
   $$
  \align
  \tilde L(s, [\frak p_2], \chi)
  &= \sum_{\frak a \in [\frak p_2], \hbox{ integral}}
     (\chi(\frak a) + \chi'(\frak a)) (N\frak a)^{-s}
     \\
    &=W(\chi) 2^{s-1}
      \sum_{b \, \odd, a \in \Bbb Z}
       \left(\frac{-D}{2a+b}\right) (2a+b -b \sqrt{-d})
        |2a+b -b \sqrt{-d}|^{-2s} .
        \endalign
        $$
 Substituting $n =2a +b$ and $m=n$,
  one gets
  $$
  \align
  \tilde L(s, [\frak p_2], \chi)
   &=W(\chi) 2^{s-1} \sum_{m, n\,  \odd} \left(\frac{-D}n\right)
         (n-m\sqrt{-d}) |n-m \sqrt{-d}|^{-2s}
         \\
         &= W(\chi) 2^s G_{\odd} (\frac{i}{2\sqrt D}, s).
         \endalign
        $$
\enddemo

{\bf Proof of Theorem 6.1}. Now it is easy to verify Theorem 6.1.
  Indeed, Lemma 6.3 gives
$$
\align &\tilde L(1, [\O_K], \chi) + \tilde L(1, [\frak p_2], \chi)
\\
 &=4 G(\frac{2i}{\sqrt D}, 1) -2 G(\frac{i}{\sqrt D}, 1)
   +2 W(\chi) G_{\odd}( \frac{i}{2\sqrt D}, 1)
   \\
    &=\frac{\pi}{\sqrt{D}}\left( 4 \theta(2) - 2 \theta(1)\right)
  +
  2W(\chi)\frac{\pi}{\sqrt{D}}\left(\theta(\frac{1}2)-\theta(1)\right)
  \\
   &=\frac{2\pi}{\sqrt{D}} (1 +W(\chi)) \left(\theta(\frac{1}2)
   -\theta(1)\right).
   \endalign
   $$
   Since $\theta(\frac{1}2) -\theta(1) >0$, Lemma 6.2 gives
   Theorem 6.1.

\subheading{7. The Central Derivative}

  Let $\chi$ be a simplest Hecke character of $K=\Bbb
  Q(\sqrt{-D})$ of conductor $\frak f$, and let
  $$
  \Lambda(s, \chi)= \left(\frac{B}{2\pi}\right)^s \Gamma(s) L(s,
  \chi),
  \tag{7.1}
  $$
  where
  $$
  B=\cases
     D &\ff  2\nmid D,
     \\
     \sqrt 2 D  &\ff 4 ||D,
     \\
      2 D &\ff 8||D.
      \endcases
      \tag{7.2}
      $$
      Then $\Lambda(s, \chi)$ has analytic continuation to the
      whole complex $s$-plane, and satisfies the functional
      equation
      $$
      \Lambda(s, \chi) = W(\chi) \Lambda(2-s, \chi).
      \tag{7.3}
      $$
      Here $W(\chi)=\pm 1$ is the root number of $\chi$. In
      section 6, we proved that $\Lambda(1, \chi) \ne 0$ if and
      only if $W(\chi)=1$. In this section, we prove

      \proclaim{Theorem 7.1} Let $\chi$ be a simplest Hecke character of $K$. Then
      $\Lambda'(1, \chi) \ne 0$ if and only if $W(\chi)=-1$.
      \endproclaim
      \demo{Proof}
             The case where $D$ is odd or divisible by $8$ was
             proved by S. Miller and the author ([MY]). We will
             settle the case $4||D$ and $D>4$ the same way. The case $D=4$ does not
             occur here. One direction is trivial by
             the functional equation. Now assume $W(\chi)=-1$. For simplicity,
             set
             $$
             L(s) =\tilde L(s, [\O_K], \chi) + \tilde L(s, [\frak
             p_2], \chi), \quad \Lambda(s)
             =\left(\frac{B}{2\pi}\right)^s \Gamma(s) L(s).
             \tag{7.4}
             $$
             Then
             $$\Lambda(s) =W(\chi)
             \Lambda(2-s)=-\Lambda(2-s).
             \tag{7.5}
             $$
             The same argument as in
             the proof of [MY, Lemma 2.1] gives
    \enddemo
             \proclaim{Lemma 7.2} When $W(\chi)=-1$ and
             $\Lambda'(1)\ne 0$, one has $\Lambda'(1, \chi)\ne 0$.
             \endproclaim

           So to prove the theorem, it suffices to show $\Lambda'(1) \ne 0$ when $W(\chi)=-1$. By Cauchy's theorem and $(7.5)$, one has
           $$
           \align
           \Lambda'(1) &=\frac{1}{2\pi i}
           \int_{2-i\infty}^{2+i\infty} \Lambda(s)
           \frac{ds}{(s-1)^2}
            -\frac{1}{2\pi i}
           \int_{-i\infty}^{i\infty} \Lambda(s)
           \frac{ds}{(s-1)^2}
       \\
        &= \frac{2}{2\pi i}\int_{2-i\infty}^{2+i\infty} \Lambda(s)
           \frac{ds}{(s-1)^2}.
           \tag{7.6}
           \endalign
           $$
  The same calculation as in the proof of Lemma 6.3 gives
  $$
  \align
  \tilde L(s, [\O_K], \chi)
   &=\sum \Sb v\, \hbox{even }, u\, \odd  \endSb \epsilon( u+ v \sqrt{-d}) (u+ v \sqrt{-d})
   |u^2+v^2 d|^{-s}
   \\
   &=2 \sum_{u>0, \, \odd}\left(\frac{-D}u\right) u^{1-2s}
     +2\sum \Sb v>0,\hbox{even}\\ u>0,\, \odd \endSb \epsilon( u+ v \sqrt{-d}) 2u
   |u^2+v^2 d|^{-s}
   \endalign
   $$
   So
   $$
\frac{1}2 \tilde L(s, [\O_K], \chi)
=\sum_{n>0}\left(\frac{-D}n\right) n^{1-2s}
  +\sum_{n>0}a(n) n^{-s}
  \tag{7.7}
  $$
  with
  $$
  a(n)=\sum \Sb u^2+v^2 d =n\\ u>0, \, \odd, v>0, \, \hbox{even}
  \endSb
        \epsilon(u+v\sqrt{-d}) 2u.
        \tag{7.8}
        $$
Similarly, Lemma 6.3(2)  gives
$$
\frac{1}2 \tilde L(s, [\frak p_2], \chi) =-\sum_{n>0} b(n)
(n/2)^{-s}, \tag{7.9}
$$
with
$$
b(n) =\sum \Sb u^2 +v^2 d=n\\ u >0, v>0, \, \odd \endSb
\left(\frac{-D}u\right) 2 u. \tag{7.10}
$$
Following [MY, pages 265-266], let
$$
f(x) =\frac{\Gamma(0, x)}x =\frac{1}x \int_x^\infty e^{-t}
\frac{dt}t \tag{7.11}
$$
be the inverse Mellin  transform of $\frac{\Gamma(s)}{(s-1)^2}$,
then
$$
\frac{1}4 \Lambda'(1) =R+C_1-C_2 \tag{7.12}
$$
where
$$
\align
 R&=\sum_{n>0} \left(\frac{-D}n\right) n f(\frac{2 \pi
n^2}B),
\\
 C_1 &= \sum_{n>0} a(n)f(\frac{2 \pi
n}B),
\\
 C_2 &= \sum_{n>0} b(n) f(\frac{ \pi
n}B).
\endalign
$$
[MY, (2.12)] asserts that the main term is
$$
R > .5235 B - .8458 B^{\frac{3}4} -.3951 B^{\frac{1}2}. \tag{7.13}
$$
The same argument as in [MY, Section 4] gives the trivial upper
bounds for $C_1$ and $C_2$ as follows:
$$
\align
 |C_1| &\le \frac{1}{4 \pi^2} \sum_{v>0} v^{-4} e^{-2 \sqrt 2 \pi
 v^2} \cdot \sum_{u>0} u e^{-\frac{2 \pi}B u^2}
 \\
  &\le  \frac{1}{4 \pi^2} \times 0.0001383 \times \frac{B}{4\pi}
  \approx 2.789\times
   10^{-7} B
  \tag{7.14}
  \endalign
$$
and
$$
\align |C_2| &\le \frac{16}{\pi^2} \sum_{v >0,\, \odd } v^{-4}
e^{-\frac{\pi v^2}{2\sqrt 2}} \sum_{u>0} u e^{-\frac{\pi u^2}B}
\\
 &\le \frac{16}{\pi^2} \times 0.3293220848 \frac{B}{2 \pi}
   \approx .0850B.
   \tag{7.15}
   \endalign
   $$
   Combining $(7.12)-7.15)$, one has
   $$
   \frac{1}4 \Lambda'(1) \ge .4385 B - .8458 B^{\frac{3}4} -.3951 B^{\frac{1}2}
   >0
   $$
when $B=\sqrt 2 D \ge 42 \sqrt 2$. This leaves one exception
$D=20$ ($D=4$ does not occur here). A simple numerical calculation
shows $\Lambda'(1)>0$ in this case too.  Now Theorem 7.1 follows
from Lemma 7.2.

 We remark that a little more work as in [MY] and [LX] would show that
 the beautiful formula $(0.5)$ also holds for
  `small' quadratic twists of the simplest Hecke characters. Similar
  results  should also hold for odd powers of these characters.
  %It
  %is interesting and important to have an explicit formula for the
  %central value and derivative of these Hecke L-functions, which
  %we do not touch in this paper. Tal Sutton is working on this.
 % Sutton is working on an explicit formula for the central
 % L-function of simplest Hecke characters in terms of theta
 % functions.

\enddocument